\documentclass[reqno]{amsart}
\usepackage{amsmath}
\usepackage{amsfonts,amssymb}
\usepackage[dvips]{graphicx}
\usepackage{epsfig}
\usepackage{color}
\newtheorem{theorem}{Theorem}[section]
\newtheorem{remark}{Remark}[section]
\newtheorem{lemma}[theorem]{Lemma}
\newtheorem{definition}{Definition}[section]

\newtheorem{proposition}[theorem]{Proposition}

\numberwithin{equation}{section}
\begin{document}

\title[Global attractor of the compressible Euler equations]
	{Existence of a global attractor for the compressible Euler equation  in a bounded interval}

\author{Yun-guang Lu}
\address{School of Mathematics, Hangzhou Normal University, Hangzhou 311121, China.}
\email{yglu2010@hznu.edu.cn}
\thanks{
Y.-G. Lu Lu's research is partially supported by the NSFC grant No.12071106 of China.
}

\author{Okihiro Sawada}
\address{Faculty of
	Engineering, Kitami Institute of Technology, 165 Koen-cho Kitami, Japan.}
\email{o-sawada@mail.kitami-it.ac.jp}

\author{Naoki Tsuge}
\address{Department of Mathematics Education, 
Faculty of Education, Gifu University, 1-1 Yanagido, Gifu
Gifu 501-1193 Japan.}
\email{tsuge.naoki.c9@f.gifu-u.ac.jp}
\thanks{
N. Tsuge's research is partially supported by Grant-in-Aid for Scientific 
Research (C) 17K05315, Japan.
}
\date{}

\maketitle

\begin{abstract}In this paper, we are concerned with the one-dimensional initial 
boundary value problem for isentropic gas dynamics.

Through the contribution of great researchers such as Lax, P. D., Glimm, J., DiPerna, R. J. and Liu, T. P., the decay theory of solutions was established. They treated with the Cauchy problem and the corresponding initial 
data have the small total variation. On the other hand,  the decay 
for initial data with large oscillation has been open for half a century. In addition, due to 
the reflection of shock waves at the boundaries, little is known for the decay of the boundary value problem on a bounded interval.

Our goal is to prove the existence of a global attractor, which 
yields a decay of solutions for large data. 
To construct approximate solutions, we introduce a modified Godunov scheme.
\end{abstract}


\section{Introduction}
The present paper is concerned with isentropic gas 
dynamics in a bounded interval $I=[0,1]$. 
\begin{align}
\begin{cases}
\displaystyle{\rho_t+m_x=0,}\\
\displaystyle{m_t+\left(\frac{m^2}{\rho}+p(\rho)\right)_x
	=0,}
\end{cases}x\in(0,1),\quad t\in(0,\infty),
\label{Euler}
\end{align}
where $\rho$, $m$ and $p$ are the density, the momentum and the 
pressure of the gas, respectively. If $\rho>0$, 
$v=m/\rho$ represents the velocity of the gas. For a barotropic gas, 
$p(\rho)=\rho^\gamma/\gamma$, where $\gamma\in(1,5/3]$ is the 
adiabatic exponent for usual gases.

We consider the initial boundary value problem (\ref{Euler}) 
with the initial and boundary data
\begin{align}  
(\rho,m)|_{t=0}=(\rho_0(x),m_0(x))\quad m|_{x=0}=m|_{x=1}=0.
\label{I.C.}
\end{align}
The above problem \eqref{Euler}--\eqref{I.C.} can be written in the following form 
\begin{align}\left\{\begin{array}{lll}
u_t+f(u)_x=0,\quad{x}\in(0,1),\quad t\in(0,\infty),\\
u|_{t=0}=u_0(x),\\
m|_{x=0}=m|_{x=1}=0
\label{IP}
\end{array}\right.
\end{align}
by using  $u={}^t(\rho,m)$, $\displaystyle f(u)={}^t\!\left(m, \frac{m^2}{\rho}+p(\rho)\right)$.

We survey the known results related to the above problem. The existence of solutions to conservation laws including \eqref{Euler} was first  established by Glimm \cite{G}. Glimm studied the Cauchy problem 
with initial data having small total variation.  {DiPerna} \cite{D2}
 proved the global existence of solutions to   
 \eqref{Euler} by the vanishing viscosity method and a compensated compactness argument. We notice that this result can treat with the arbitrary $L^{\infty}$ data.

The theory of decay for genuinely nonlinear 2$\times$2 systems of conservation laws was established by Glimm-Lax \cite{GL}. Glimm-Lax showed that if initial data are constant outside a finite interval and have locally bounded total variation and small oscillation, then the tonal 
variation of the solution of \cite{G} decays to zero. 
The Glimm-Lax theory had been further developed by DiPerna and Liu: 
general conservation laws with a convex entropy function \cite{D1}, 
general conservation laws with small initial data in total variation \cite{L}. However, the decay for initial data with large oscillation has been open for half a century. Recently, Tsuge \cite{T10} 
obtained a decay estimate for large $L^{\infty}$ data. 
On the other hand, little is known for the decay of the 
boundary value problem such as \eqref{IP}, because it is difficult 
to treat with the reflection of shock waves at the both boundaries.

Our goal in the present paper is to investigate the decay structure of the boundary value problem \eqref{IP}. We first introduce a modified Godunov scheme developed in \cite{T1}--\cite{T10}
to construct approximate solutions. To deduce the convergence, 
we employ the compensated compactness developed by DiPerna. We next prove the existence of a global attractor, which 
yields decay estimates of solutions for initial data with large oscillation.

To state our main theorem, we define the Riemann invariants $w,z$, which play important roles
in this paper, as
\begin{definition}
	\begin{align*}
	w:=\frac{m}{\rho}+\frac{\rho^{\theta}}{\theta}=v+\frac{\rho^{\theta}}{\theta},
	\quad{z}:=\frac{m}{\rho}-\frac{\rho^{\theta}}{\theta}
	=v-\frac{\rho^{\theta}}{\theta}\quad
	\left(\theta=\frac{\gamma-1}{2}\right).
	\end{align*}
\end{definition}
These Riemann invariants satisfy the following.
\begin{remark}\label{rem:Riemann-invariant}
	\normalfont
	\begin{align*}
	&|w|\geq|z|,\;w\geq0,\;\mbox{\rm when}\;v\geq0.\quad
	|w|\leq|z|,\;z\leq0,\;\mbox{\rm when}\;v\leq0.
    \\
	&v=\frac{w+z}2,
	\;\rho=\left(\frac{\theta(w-z)}2\right)^{1/\theta},\;m=\rho v.
	\end{align*}From the above, the lower bound of $z$ and the upper bound of $w$ yield the bound of $\rho$ and $|v|$.
\end{remark}

Moreover, we define the entropy weak solution.
\begin{definition}
		A measurable function $u(x,t)$ is called an {\it entropy weak solution} of the initial boundary value problem \eqref{IP} if 
	\begin{align*}
		&\int^{1}_{0}\int^{\infty}_0\rho(\varphi_1)_t+m(\varphi_1)_xdxdt+\int^{1}_{0}
		\rho_0(x)\varphi_1(x,0)dx=0,\\
&\int^{1}_{0}\int^{\infty}_0m(\varphi_2)_t+\left(\frac{m^2}{\rho}+p(\rho)\right)(\varphi_2)_x dxdt+\int^{1}_{0}
		m_0(x)\varphi_2(x,0)dx=0
	\end{align*}
	hold for any test function $\varphi_1,\varphi_2\in C^1_0([0,1]\times[0,\infty))$ 
satisfying $\varphi_2(0,t)=\varphi_2(1,t)=0$ and 
	\begin{align}
		&\int^{1}_0\int^{\infty}_0\hspace{-1ex}\eta(u)\psi_t+q(u)\psi_xdxdt\geq0
		\label{entropy solution}
	\end{align}
	holds for any non-negative test function $\psi\in C^1_0((0,1)\times(0,\infty))$, where 
	$(\eta,q)$ is a pair of convex entropy--entropy flux of \eqref{Euler}.
\end{definition}

We set 
$\displaystyle \bar{\rho}=\int^1_0\rho_0(x)dx
$. Since $\displaystyle \bar{\rho}=0$ implies that the initial data becomes vacuum, we assume $\displaystyle \bar{\rho}>0$. In addition, we define the mechanical energy as 
\begin{align*}
\eta_{\ast}(u)=\frac12\frac{m^2}{\rho}+\frac1{\gamma(\gamma-1)}\rho^{\gamma}.
\end{align*}
We choose a positive constant $\mu$ small enough. We then set
\begin{align}
\begin{alignedat}{2}
&\bar{\eta}=\int^1_0\eta(u_0(x))dx+\mu,\quad
\nu=\dfrac{3\gamma-1}{\gamma+1}\dfrac{\bar{\eta}}{\bar{\rho}},\quad
K=\bar{\rho}\nu-\bar{\eta},\\
&M_{\infty}=\dfrac{4}{3\gamma-1}
\left(\dfrac{2\gamma^2(\gamma-1)}{3\gamma-1}\right)^{\frac{\gamma+1}{2(\gamma-1)}}\dfrac{\nu^{\frac{3\gamma-1}{2(\gamma-1)}}}{\bar{\rho}\nu-\bar{\eta}}+\dfrac{2\left(\nu\bar{\rho}+\bar{\eta}+K\right)}{\gamma-1}.
\label{notation}
\end{alignedat}
\end{align}
We notice that $K>0$, if necessary, by choosing $\mu$ small enough.

Moreover, for any fixed positive constant $\varepsilon'$, we define 
\begin{align}
	\begin{split}
		&\tilde{z}(x,t)=z(x,t)-\varepsilon'\int^x_{0}\zeta(u(y,t))dy,
\;
\tilde{w}(x,t)=w(x,t)-\varepsilon'\int^x_{0}\zeta(u(y,t))dy,
	\end{split}
	\label{transformation}
\end{align}
where 
\begin{align}
\zeta(u)=\eta_{\ast}(u)-\nu\rho+K.
\label{zeta}
\end{align}
From the conservation of mass and the energy inequality, we find that 
\begin{align*}
\int^x_{0}\zeta(u(y,t))dy\leq& 
\int^x_0 \left\{\eta_{\ast}(u(y,t))+\nu\rho(y,t)+K\right\}dy\\
\leq& \int^1_0 \left\{\eta_{\ast}(u(x,t))+\nu\rho(x,t)+K\right\}dy
\\
\leq& \int^1_0 \left\{\eta_{\ast}(u_0(x))+\nu\rho_0(x)+K\right\}dy.
\end{align*}

Our main theorem is as follows.
\begin{theorem}\label{thm:main}We assume that 
\begin{align}
\begin{split}
\rho_0(x)\geq0\quad{a.e.\ }x\in{I},\quad \rho_0\in L^{\infty}(I),\quad\dfrac{m_0}{\rho_0}\in L^{\infty}(I).
\end{split}
\label{maintheorem1}
\end{align}

Then, there exists a global entropy weak solution of the initial boundary value problem \eqref{IP}. Moreover, for any positive constant $\varepsilon$, 
there exists positive constant $t_0$ such that the solution satisfies 
\begin{align}
\begin{split}
&-M_{\infty}-\varepsilon\leq \tilde{z}(x,t),\quad 
\tilde{w}(x,t)\leq M_{\infty}+\varepsilon,\quad \rho(x,t)\geq0,\\
&{a.e.\ }(x,t)\in{I}\times[t_0,\infty),
\end{split}
\label{maintheorem2}
\end{align}where $t_0$ depends only on $\varepsilon,\;\varepsilon'$ and the bound of initial data.

\end{theorem}	
For simplicity, we set $\varepsilon'=1$ hereafter.

\begin{remark}\label{rem:boundary}
In this remark, we state the important conditions necessary to construct an invariant region at the boundary. This condition 
will be used in Section 3 \eqref{I_2N}.

We choose a positive constant $M_0$ such that 
\begin{align}
-M_0+\int^x_{0}\zeta(u_0(y))dy\leq z(u_0(x)),\;
w(u_0(x))\leq M_0+\int^x_{0}\zeta(u_0(y))dy.\label{inequaltiy1}
\end{align}
Then, in the proof of Theorem \ref{thm:main}, we will observe that there exists a continuous function 
${\mathcal M}(t)$ such that ${\mathcal M}(0)=M_0$,
${\mathcal M}(t_0)=M_{\infty}$ and 
\begin{align}
-{\mathcal M}(t)+\int^x_{0}\zeta(u(y,t))dy\leq z(u(x,t)),\;
w(u(x,t))\leq {\mathcal M}(t)+\int^x_{0}\zeta(u(y,t))dy.
\label{inequaltiy2}
\end{align}
Let the lower and upper bounds in \eqref{inequaltiy2} be 
	\begin{align*}
			L(x,t;u)=-{\mathcal M}(t)+\int^x_{0}\zeta(u(y,t))dy,\;
			U(x,t;u)={\mathcal M}(t)+\int^x_{0}\zeta(u(y,t))dy,
	\end{align*}
respectively. Then we notice that 
	\begin{align}
		-L(0,t;u)\leq U(0,t;u),\;-L(1,t;u)\geq U(1,t;u).
		\label{boundary-bound}
	\end{align}
In fact, the former is clear. The latter is from \eqref{notation} and the energy inequality deduced as follows.
	\begin{align}
		L(1,t;u)+U(1,t;u)=&2\int^1_{0}\left\{\eta_{\ast}(u(x,t))-\nu\rho(x,t)+K\right\}dx\nonumber\\=&2\int^1_{0}\left\{\eta_{\ast}(u(x,t))-\nu\rho(x,t)+\nu\bar{\rho}-\bar{\eta}\right\}dx\nonumber\\
		=&2\int^1_{0}\left\{\left(\eta_{\ast}(u(x,t))-\eta_{\ast}({u}_0(x))\right)
		-\nu\left(\rho(x,t)-\bar{\rho}\right)-\mu\right\}dx\nonumber\\
		\leq&-2\mu.
\label{boundary-bound2}
	\end{align}

\eqref{boundary-bound} is a necessary condition that 
\eqref{inequaltiy2} holds for boundary data $m=0$. 
\end{remark}

\subsection{Outline of the proof (formal argument)}

The proof of the main theorem is a little complicated. Therefore, 
before proceeding to the subject, let us grasp the point of the main estimate by a formal argument. 
We assume that a solution is smooth and the density is nonnegative in this section.

We consider the physical region $\rho\geq0$ (i.e., $w\geq z$.). Recalling Remark \ref{rem:Riemann-invariant}, it suffices to 
derive the lower bound of $z(u)$ and the upper bound of $w(u)$ to obtain the bound of $u$. To do this, we diagonalize \eqref{Euler}. 
If solutions are smooth, we deduce from \eqref{Euler} 
\begin{align}
z_t+\lambda_1z_x=0,\quad
w_t+\lambda_2w_x=0,
\label{Euler2}
\end{align} 
where $\lambda_1$ and $\lambda_2$ are the characteristic speeds defined as follows 
\begin{align}
\lambda_1=v-\rho^{\theta},\quad\lambda_2=v+\rho^{\theta}.
\label{char}
\end{align}

We introduce $\tilde{z},\tilde{w},\tilde{\rho},\tilde{v},\tilde{\lambda}_1,\tilde{\lambda}_2
$ as follows.
\begin{align}
	\begin{split}
		&z=\tilde{z}+\int^x_{0}\left\{\eta_{\ast}(u)-\nu\rho+K\right\}dy,\quad w=\tilde{w}+\int^x_{0}\left\{\eta_{\ast}(u)-\nu\rho+K\right\}dy,\\
        &\tilde{\rho}=\left(\dfrac{\theta(\tilde{w}-\tilde{z})}{2}\right)^{1/\theta},\quad \tilde{v}=\dfrac{\tilde{w}+\tilde{z}}{2},\quad
\tilde{\lambda}_1=\tilde{v}-\tilde{\rho}^{\theta},\quad
\tilde{\lambda}_2=\tilde{v}+\tilde{\rho}^{\theta}.
	\end{split}
	\label{transformation}
\end{align}

We denote the flux of $\eta_{\ast}(u)$ by
\begin{align}
q_{\ast}(u)=m\left(\frac12\frac{m^2}{\rho^2}+\frac{\rho^{\gamma-1}}{\gamma-1}\right). 
\label{energy-flux}
\end{align}
Then, from \eqref{Euler}, it holds that 
\begin{align}
\left(\eta_{\ast}(u)\right)_t+\left(q_{\ast}(u)\right)_x=0.
\label{energy-conservation}
\end{align}

For $\delta=\theta K\varepsilon/2$, we define 
$\hat{z}=\tilde{z}-\delta t,\;\hat{w}=\tilde{w}+\delta t$.
We then deduce from $\eqref{Euler}_1$ and \eqref{energy-conservation} that
\begin{align}
	\hat{z}_t+\lambda_1 \hat{z}_x=g_1(x,t,u),\quad
	\hat{w}_t+\lambda_2 \hat{w}_x=g_2(x,t,u),
	\label{Riemann2}
\end{align}
where 
\begin{align}
	\begin{alignedat}{2}
		&g_1(x,t,u)=&&-K\lambda_1+\dfrac{1}{\gamma(\gamma-1)}\rho^{\gamma+\theta}
		+\dfrac{1}{\gamma}\rho^{\gamma}v+\dfrac{1}{2}\rho^{\theta+1}v^2-\nu\rho^{\theta+1}-\delta,\\
		&g_2(x,t,u)=&&-K\lambda_2-\dfrac{1}{\gamma(\gamma-1)}\rho^{\gamma+\theta}
		+\dfrac{1}{\gamma}\rho^{\gamma}v-\dfrac{1}{2}\rho^{\theta+1}v^2+\nu\rho^{\theta+1}+\delta.
	\end{alignedat}
\label{inhomo2}
\end{align}
On the other hand, we notice that
\begin{align*}
	-M_0\leq \hat{z}_0(x),\; \hat{w}_0(x)\leq M_0.
\end{align*}
Our goal is to prove that  
\begin{align}
	\hat{S}_{inv}=\{(\hat{z},\hat{w})\in{\bf R}^2;-M_0\leq\hat{z},\;\hat{w}\leq M_0\}
\label{invariant region1}
\end{align}
is an invariant region for $0\leq t\leq t_0$, where $t_0=
\max\{(M_0-M_{\infty}-\varepsilon)/\delta,0\}$. 

We consider the case where $M_0>M_{\infty}+\varepsilon$. To achieve this, assuming that 
\begin{align*}
	-M_0< \hat{z}_0(x),\; \hat{w}_0(x)< M_0
\end{align*}
and there exist $x_{\ast}\in(0,1),\;0<t_{\ast}\leq t_0$ such that 
\eqref{invariant1} or \eqref{invariant2} holds, we shall
deduce a contradiction, where
\begin{align}
\begin{alignedat}{2}
	&-M_0<\hat{z}(x,t),\;\hat{w}(x,t)< M_0,\quad x\in(0,1),\;
	0\leq t<t_{\ast}\\&\text{\hspace*{0ex}and}\quad
	\hat{z}(x_{\ast},t_{\ast})=-M_0,\;\hat{w}(x_{\ast},t_{\ast})\leq M_0,	
\end{alignedat}\label{invariant1}\\	
\begin{alignedat}{2}
	&-M_0<\hat{z}(x,t),\;\hat{w}(x,t)< M_0,\quad x\in(0,1),\;
	0\leq t<t_{\ast}\\&\text{\hspace*{0ex}and}\quad
	-M_0\leq\hat{z}(x_{\ast},t_{\ast}),\;\hat{w}(x_{\ast},t_{\ast})=M_0.
\end{alignedat}\label{invariant2}
\end{align}To do this, we prove 
\begin{align}	&g_1(x_{\ast},t_{\ast},u)>0\text{, when \eqref{invariant1} holds},
	\label{g1}\\
	&g_2(x_{\ast},t_{\ast},u)<0\text{, when \eqref{invariant2} holds}.
	\label{g2}
\end{align}

From \eqref{transformation}, we notice $\tilde{\rho}=\rho$. We thus obtain
\begin{align}
	\dfrac{(\rho(x,t))^{\theta}}{\theta}=&\dfrac{({\tilde{\rho}(x,t)})^{\theta}}{\theta}
=\dfrac{\tilde{w}(x,t)-\tilde{z}(x,t)}2
=\dfrac{\hat{w}(x,t)-\hat{z}(x,t)-2\delta t}2
\leq M_0-\delta t
\label{eqn:upper_rho}
\end{align}
and observe 
\begin{align}
	\begin{alignedat}{2}&\lambda_1=z+\dfrac{3-\gamma}{\gamma-1}\rho^{\theta}
=\tilde{z}+\int^x_0\zeta(u)dx+\dfrac{3-\gamma}{\gamma-1}\rho^{\theta},\\
		&\lambda_2=w-\dfrac{3-\gamma}{\gamma-1}\rho^{\theta}
=\tilde{w}+\int^x_0\zeta(u)dx-\dfrac{3-\gamma}{\gamma-1}\rho^{\theta}.
\label{eqn:lambda_2}
	\end{alignedat}
\end{align}

For $(x,t)=(x_{\ast},t_{\ast})$, since  
$M_0-\delta t_{\ast}\geq M_0-\delta t_0=M_{\infty}
+\varepsilon$, recalling $\delta=\theta K\varepsilon/2$ and \eqref{notation}, we deduce from $\eqref{eqn:upper_rho}$ and $\eqref{eqn:lambda_2}$ that  
	\begin{align}
		\begin{alignedat}{3}
			&g_1(x,t,u)&&=&&-K\left(\tilde{z}+\int^x_0\zeta(u)dx
+\dfrac{3-\gamma}{\gamma-1}\rho^{\theta}\right)
+\dfrac{\rho^{\theta+1}}2\left({v}+\frac{\rho^{\theta}}{\gamma}\right)^2
			+\dfrac{\gamma+1}{2\gamma^2(\gamma-1)}\rho^{\gamma+\theta}\\
&&&&&-\nu\rho^{\theta+1}-\delta\\
&&&\geq&&K\left(M_0-\delta t_{\ast}\right)-K\left(\nu\bar{\rho}+\bar{\eta}+K\right)
-\dfrac{3-\gamma}{\gamma-1}\theta K\left(M_0-\delta t_{\ast}\right)\\
&&&&&+\min_{\rho}\left\{\dfrac{\gamma+1}{2\gamma^2(\gamma-1)}\rho^{\gamma+\theta}-\nu\rho^{\theta+1}
			\right\}-\delta\\
&&&\geq&&\theta K\left(M_{\infty}
+\varepsilon\right)-K\left(\nu\bar{\rho}+\bar{\eta}+K\right)\\&&&&&
+\dfrac{2(\gamma-1)}{3\gamma-1}
\left(\dfrac{2\gamma^2(\gamma-1)}{3\gamma-1}\right)^{\frac{\gamma+1}{2(\gamma-1)}}\nu^{\frac{3\gamma-1}{2(\gamma-1)}}-\delta\\
&&&=&&\delta\\
&&&>&&0.
			\end{alignedat}
		\label{estimate1}
	\end{align}
On the other hand, since $\tilde{z}$ attains the minimum at $(x,t)=(x_{\ast},t_{\ast})$, we find that 
$\tilde{z}_t(x_{\ast},t_{\ast})\leq0,\;\tilde{z}_x(x_{\ast},t_{\ast})=0$. Then, from $\eqref{Riemann2}_1$, 
we have $g_1(x,t,u)>0$ at $(x,t)=(x_{\ast},t_{\ast})$. This is a contradiction for \eqref{estimate1}.

	For $(x,t)=(x_{\ast},t_{\ast})$, we can similarly obtain  
	\begin{align}
		\begin{alignedat}{3}
			&g_2(x,t,u)&&=&&-K\left(\tilde{w}+\int^x_0\zeta(u)dx-\dfrac{3-\gamma}{\gamma-1}\rho^{\theta}\right)-\dfrac{\rho^{\theta+1}}2\left({v}-\frac{\rho^{\theta}}{\gamma}\right)^2
			-\dfrac{\gamma+1}{2\gamma^2(\gamma-1)}\rho^{\gamma+\theta}\\
&&&&&+\nu\rho^{\theta+1}+\delta\\
&&&\leq&&-K\left(M_0-\delta t_{\ast}\right)+K\left(\nu\bar{\rho}+\bar{\eta}+K\right)
+\dfrac{3-\gamma}{\gamma-1}\theta K\left(M_0-\delta t_{\ast}\right)\\
&&&&&+\max_{\rho}\left\{-\dfrac{\gamma+1}{2\gamma^2(\gamma-1)}\rho^{\gamma+\theta}+\nu\rho^{\theta+1}
			\right\}+\delta\\
&&&\leq&&-\theta K\left(M_{\infty}
+\varepsilon\right)+K\left(\nu\bar{\rho}+\bar{\eta}+K\right)\\&&&&&
+\dfrac{2(\gamma-1)}{3\gamma-1}
\left(\dfrac{2\gamma^2(\gamma-1)}{3\gamma-1}\right)^{\frac{\gamma+1}{2(\gamma-1)}}\nu^{\frac{3\gamma-1}{2(\gamma-1)}}+\delta\\
&&&=&&-\delta\\
&&&<&&0.
			\end{alignedat}
		\label{estimate2}
	\end{align}
Since $\tilde{w}$ attains the maximum at $(x,t)=(x_{\ast},t_{\ast})$, we find that 
$\tilde{w}_t(x_{\ast},t_{\ast})\geq0,\;\tilde{w}_x(x_{\ast},t_{\ast})=0$. Then, from $\eqref{Riemann2}_2$, 
we have $g_2(x,t,u)>0$ at $(x,t)=(x_{\ast},t_{\ast})$. This is a contradiction.

\eqref{invariant region1} implies that $(\tilde{z}(x,t_0),
\tilde{w}(x,t_0))$ is contained in 
\begin{align*}
	\tilde{S}_{inv}=\{(\tilde{z},\tilde{w})\in{\bf R}^2;-M_{\infty}-\varepsilon\leq\tilde{z},\;\tilde{w}
\leq M_{\infty}+\varepsilon\}.
\end{align*}
In addition, we find that $\tilde{S}_{inv}$ is an invariant region in the similar manner to \eqref{invariant region1}. Therefore, we 
can prove \eqref{maintheorem2}.

The present paper is organized as follows.
In Section 2, we construct approximate solutions by 
the Godunov scheme mentioned above. In Section 3, we drive the bounded a decay estimate of our approximate solutions.

\section{Construction of Approximate Solutions}
\label{sec:construction-approximate-solutions}
In this section, we construct approximate solutions. In the strip 
$0\leq{t}\leq [\hspace{-1.2pt}[T]\hspace{-1pt}]+1$ for any fixed positive constant $T$, we denote these 
approximate solutions by $u^{\varDelta}(x,t)
=(\rho^{\varDelta}(x,t),m^{\varDelta}(x,t))$, where $[\hspace{-1.2pt}[T]\hspace{-1pt}]$ is the greatest integer 
not greater than $T$.  
For $N_x\in{\bf N}$, we define the space mesh
lengths by ${\varDelta}x=1/(2N_x)$. We take time mesh length ${\varDelta}{t}$ such that 
\begin{align}
	\frac{{\varDelta}x}{{\varDelta}{t}}=2[\hspace{-1.2pt}
[\max\{M_0,M_{\infty}+\varepsilon\}+\bar{\eta}+\nu\bar{\rho}+K]\hspace{-1pt}]+1,
\label{CFL}
\end{align}
where $[\hspace{-1.2pt}[x]\hspace{-1pt}]$ is the greatest integer 
not greater than $x$. Then we define $N_t=([\hspace{-1.2pt}[T]\hspace{-1pt}]+1)/(2{\varDelta t})\in{\bf N}$.
In addition, 
we set 
\begin{align*}
	(j,n)\in{\bf N}_x\times{\bf N}_t,
\end{align*}
where ${\bf N}_x=\{1,3,5,\ldots,2N_x-1\}$ and ${\bf N}_t=\{0,1,2,\ldots,2N_t\}$.  For simplicity, we use the following terminology
\begin{align}
\begin{alignedat}{2}
&x_j=j{\varDelta}x,\;t_n=n{\varDelta}t,\;t_{n.5}=\left(n+\frac12\right){\varDelta}t,
\;t_{n-}=n{\varDelta}t-0,\;t_{n+}=n{\varDelta}t+0.
\end{alignedat}
\label{terminology} 
\end{align}

First we define $u^{\varDelta}(x,-0)$ by $u^{\varDelta}(x,-0)=u_0(x)$. Then, for $j\in {\bf N}_x$, we denote $E_j^0(u)$ by
\begin{align*}
	E_j^0(u)=\frac1{2{\varDelta}x}\int^{x_{j+1}}_{x_{j-1}}
	u^{\varDelta}(x,-0)dx.
\end{align*}

Next, assume that $u^{\varDelta}(x,t)$ is defined for $t<{t}_{n}$. 
Then, for $j\in {\bf N}_x$, we denote $E^n_j(u)$ by 
\begin{align*}
	E^n_j(u)=\frac1{2{\varDelta}x}\int_{x_{j-1}}^{x_{j+1}}u^{\varDelta}(x,t_{n-})dx.
\end{align*}
Let $E^n(x;u)$ be a piecewise constant function defined by
\begin{align*}
	E^n(x;u)=
		E^n_j(u),\quad &x\in [x_{j-1},x_{j+1})\quad (j\in {\bf N}_x).
\end{align*}

To define $u_j^n=(\rho_j^n,m_j^n)$ for $j\in {\bf N}_x$, we first determine symbols $I^n_j$ and $L_n$. Let the approximation of $\zeta(u)$ be
\begin{align*}
	I^n_j:=\int^{x_{j-1}}_{0}\zeta(E^n(x;u))dx 
		+\frac{1}2\int_{x_{j-1}}^{x_{j+1}}\zeta(E^n(x;u))dx
         =\int^{x_{j}}_{0}\zeta(E^n(x;u))dx,
\end{align*}where $\zeta$ is defined in \eqref{zeta}.

Let ${\mathcal D}=(x(t),t)$ denote a discontinuity in 
${u}^{\varDelta}(x,t),\;[\eta_{\ast}]$ and $[q_{\ast}]$ 
denote the jump of $\eta_{\ast}({u}^{ \varDelta}(x,t))$ and $q_{\ast}({u}^{ \varDelta}(x,t))$ across ${\mathcal D}$ from 
left to right, respectively,
\begin{align*}
&[\eta_{\ast}]=\eta_{\ast}({u}^{\varDelta}(x(t)+0,t))-\eta_{\ast}({u}^{ \varDelta}(x(t)-0,t)),
\\&
{[q_{\ast}]=q_{\ast}({u}^{ \varDelta}(x(t)+0,t))-q_{\ast}({u}^{ \varDelta}(x(t)-0,t))},
\end{align*}
where $q_{\ast}(u)$ is defined in \eqref{energy-flux}.

To measure the error in the entropy condition and the gap of the 
energy at $t_{n\pm}$, we introduce the following functional.
\begin{align}
	\begin{alignedat}{2}
	L_n=&\int^{t_{n}}_{0}\sum_{0\leq x\leq 1}\left(\sigma[\eta_{\ast}]-[q_{\ast}]\right)dt
	+\sum^n_{k=0}\int^1_0\left\{\eta_{\ast}({u}^{ \varDelta}(x,t_{k-0}))-\eta_{\ast}(E^k(x;u))\right\}dx\\&
	+\sum^n_{k=0}\sum_{j\in J_k}\frac1{2{\varDelta}x}\int^{x_{j+1}}_{x_{j-1}}\int^{x}_{x_{j-1}}
	R^k_j(y)dydx,
	\end{alignedat}
\label{functional discontinuity}	
\end{align}
where
\begin{align*}
R^n_j(x)=&\int^1_0(1-\tau)\cdot{}^t\left({u}^{\varDelta}(x,t_{n-})-E^n(x;u)\right)\\&\times
\nabla^2\eta_{\ast}\left(E^n(x;u)+\tau\left\{{u}^{\varDelta}(x,t_{n-})-E^n(x;u)\right\}\right)\left({u}^{\varDelta}(x,t_{n-})-E^n(x;u)\right)d\tau
\end{align*}
and the summention in $\sum_{0\leq x\leq1}$
is taken over all discontinuities in ${u}^{ \varDelta}(x,t)$ at a fixed time $t$ over 
$x\in[0,1]$,
$\sigma$ is the propagating speed of the discontinuities.

From the entropy condition, $\sigma[\eta_{\ast}]-[q_{\ast}]\geq0$. From the Jensen inequality, 
$\int^1_0\left\{\eta_{\ast}({u}^{ \varDelta}(x,t_{n-0}))-\eta_{\ast}(E^n(x;u))\right\}dx\geq0$. Therefore, we find that $L_n\geq0$.

Using $I^n_j$ and $L_n$, we define $u_j^n$ as follows. First, we define a sequence $\left\{M_n\right\}_{n\in {\bf N}_t}$ with the initial term $M_0$ as follows.
\begin{align} 
M_{n+1}=
\begin{cases}
M_n-\delta{\varDelta}t,&\text{when }M_{n}+L_n\geq M_{\infty}+\varepsilon,\\
M_n,&\text{when }M_{n}+L_n<M_{\infty}+\varepsilon.
\end{cases}
\label{M_n}
\end{align}We notice that $M_n+L_n\geq
M_{\infty}+\varepsilon-\delta{\varDelta}t$.

Next, we choose $\beta$ such that $1<\beta<1/(2\theta)$. If 
\begin{align*}
	E^n_j(\rho):=
	\frac1{2{\varDelta}x}\int_{x_{j-1}}^{x_{j+1}}\rho^{\varDelta}(x,t_{n-})dx<({\varDelta}x)^{\beta},
\end{align*} 
we define $u_j^n$ by $u_j^n=(0,0)$;
otherwise, setting\begin{align}
		{z}_j^n:=
		\max\left\{z(E_j^n(u)),\;-M_n-L_n+I^n_j\right\},\;
		w_j^n:=\min\left\{w(E_j^n(u)),\;M_n+L_n+I^n_j\right\}
		,
	\label{def-u^n_j}
\end{align}

we define $u_j^n$ by
\begin{align*}
	u_j^n:=(\rho_j^n,m_j^n):=(\rho_j^n,\rho_j^nv^n_j)
	:=\left(\left\{\frac{\theta(w_j^n-z_j^n)}{2}\right\}
	^{1/\theta},
	\left\{\frac{\theta(w_j^n-z_j^n)}{2}\right\}^{1/\theta}
	\frac{w_j^n+z_j^n}{2}\right).
\end{align*}

\begin{remark}\label{rem:E}\normalfont
We find 
	\begin{align}
		\begin{split}
			-M_n-L_n+I^n_j\leq z(u_j^n),\quad
			{w}(u_j^n)\leq M_n+L_n+I^n_j.
		\end{split}
		\label{remark2.1}
	\end{align}

	This implies that we cut off the parts where 
	$z(E_j^n(u))<-M_n-L_n+I^n_j$
	and $w(E_j^n(u))>M_n+L_n+I^n_j$
	in  defining $z(u_j^n)$ and 
	${w}(u_j^n)$. 
\end{remark}

We must construct our approximate solutions $u^{\varDelta}(x,t)$ near the boundary and in an interior domain. Here we recall \eqref{boundary-bound}. This condition is necessary so that the inequality \eqref{inequaltiy2} holds even if a shock wave appears at the boundary $x=1$. To see this, we are devoted to treating with the construction near the boundary $x=1$. 
For the construction in the interior domain, refer to \cite{T9}.

\label{subsec:construction-approximate-solutions}
We then assume that 
approximate solutions $u^{\varDelta}(x,t)$ are defined in domains $D_1: 
t<{t}_n\quad(
n\in {\bf N}_t)$ and $D_2:
x<x_{2N_x-1},\;{t}_n\leq{t}<t_{n+1}$. 
By using $u_j^n$ defined above and $u^{\varDelta}(x,t)$ defined in $
D_2$, we 
construct the approximate solutions  in the cell $D:\;{t}_n\leq{t}<{t}_{n+1}\quad 
(n\in{\bf N}_t),\; x_{2N_x-1}\leq{x}<x_{2N_x}$.

\begin{figure}[htbp]
	\begin{center}
		\hspace{-2ex}
		\includegraphics[scale=0.4]{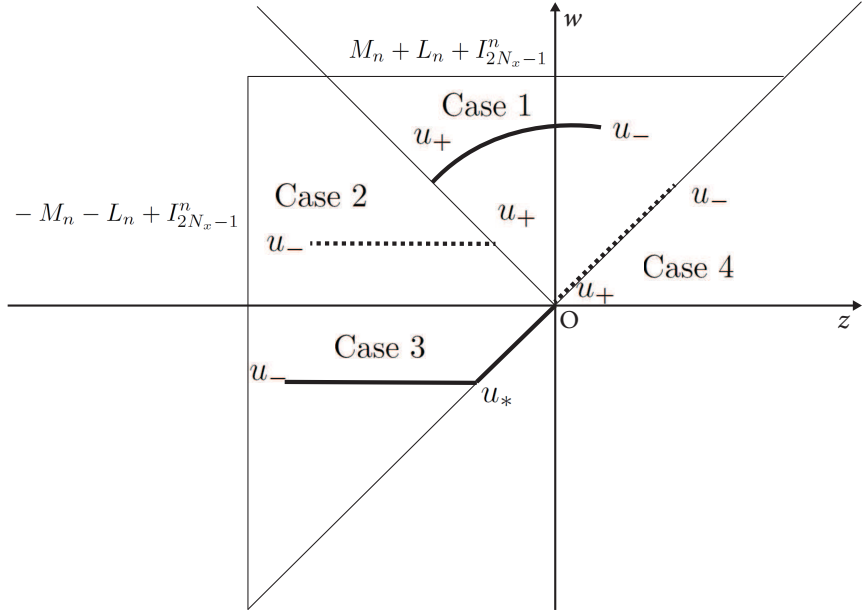}
	\end{center}
	\caption{Solutions of the Riemann initial boundary value problem \eqref{Euler} and  \eqref{eqn:RiemannIB} in $(z,w)$-plane.}
	\label{Fig:Rimann}
\end{figure}

We denote $u^n_{2N_x-1}$ by $u_-=(\rho_-,m_-)=(\rho_-,\rho_-v_-)$ and solve the Riemann initial boundary value problem \eqref{Euler} and 
\begin{equation}
u|_{t=t_n}=u_-,\quad{m}|_{x=1}=0
\label{eqn:RiemannIB}
\end{equation} 
in $D$. We draw a diagram by using 
the wave curve of the first 
family and the vacuum as follows (see \cite{T1} and Figure \ref{Fig:Rimann}): 
\begin{itemize}
\item[Case 1] If $\rho_->0$ and $v_-\geq0$, there exists $u_+=(\rho_+,m_+)=(\rho_+,\rho_+v_+)$ with $v_+=0$ from which $u_+$ is connected by a 1-shock curve.
\item[Case 2] If $v_-\leq0$ and $w(u_-)\geq0$, then there exists $u_+$ with $v_+=0$ from which $u_+$ is connected by a 1-rarefaction curve. 
\item[Case 3] If $v_-\leq0$ and $w(u_-)\leq0$, then there exists $u_*$ with $\rho_*=0$ 
from which $u_-$ is connected by a 1-rarefaction, and $u_*$ and $u_+$ with 
$\rho_+=v_+=0$ are connected 
by the vacuum.
\item[Case 4] If $v_-\geq0$ and $\rho_-=0$, then $u_+$ with   
$\rho_+=v_+=0$ is connected from $u_-$ by the vacuum.
\end{itemize}

\subsection{Case 1: the case where a shock waves arise}
In this case, 
the solution $u_{\rm R}(x,t)$ of \eqref{Euler} and \eqref{eqn:RiemannIB} is as follows. 
\begin{align*}u_{\rm R}(x,t)=
\begin{cases}
u_-,&D\cap\left\{x-x_{2N_x}\leq\sigma_s(t-t_n)\right\},\\ 
u_+,&D\cap\left\{x-x_{2N_x}>\sigma_s(t-t_n)\right\},
\end{cases}
\end{align*}
where $\sigma_s$ is the speed of 1-shock wave.

We next replace the above constant states $u_-,\;u_+$ with functions of $x$ and $t$ as follows:

In view of \eqref{transformation}, we construct ${u}^{\varDelta}_-(x,t)$. We first determine the approximation of $\tilde{z},\tilde{w}$ in \eqref{transformation} 
as follows.
\begin{align*}
	\begin{alignedat}{2}
		\tilde{z}^{\varDelta}_-
		=&z_--
\int^{x_{2N_x-1}}_{0}
		\zeta(u^{\varDelta}_{n,0}(x))dx,\;
		\tilde{w}^{\varDelta}_1=w_--
 \int^{x_{2N_x-1}}_{0}
		\zeta(u^{\varDelta}_{n,0}(x))dx,
	\end{alignedat}
	\end{align*}
where $u^{\varDelta}_{n,0}(x)$
is a piecewise constant function defined by
\begin{align}
	u^{\varDelta}_{n,0}(x)=
		u^n_j,\quad &x\in [x_{j-1},x_{j+1})\quad (j\in {\bf N}_x).	
	\label{def-u0}
\end{align} 
We set       
\begin{align}
	\begin{alignedat}{2}
		&\check{z}^{\varDelta}_-(x,t)=&&\tilde{z}^{\varDelta}_-
		+\int^{x_{2N_x-1}}_{0}
		\zeta(u^{\varDelta}_{n,0}(x))dx
        +\int^x_{x_{2N_x-1}}\zeta(u_-)dy
\\&&&	
+\left\{g_1(x,t;u_-)+V(u_-)\right\}(t-t_n)
		,\\
		&\check{w}^{\varDelta}_-(x,t)=&&\tilde{w}^{\varDelta}_-
		+ \int^{x_{2N_x-1}}_{0}
		\zeta(u^{\varDelta}_{n,0}(x))dx+\int^x_{x_{2N_x-1}}\zeta(u_-)dy
\\&&&
+\left\{g_2(x,t;u_-)+V(u_-)\right\}
		(t-t_n),
	\end{alignedat}\label{appro1-2}
\end{align}
where $g_1$ and $g_2$ are defined in \eqref{inhomo2}, 
\begin{align}
V(u)=q_{\ast}(u)-\nu m.
\label{V}
\end{align}

Using $\check{u}^{\varDelta}_-(x,t)$, we next define ${u}^{\varDelta}_-(x,t)$ as follows. 
\begin{align}
	\begin{alignedat}{2}
		&{z}^{\varDelta}_-(x,t)=&&\tilde{z}^{\varDelta}_-
+\int^{x_{2N_x-1}}_{0}
		\zeta(u^{\varDelta}_{n,0}(x))dx+\int^x_{x_{2N_x-1}}
\zeta(\check{u}^{\varDelta}_-(y,t))dy\\&&&	+
\left\{g_1(x,t;\check{u}^{\varDelta}_-)+V(u_-)\right\}(t-t_n)
		,\\
		&{w}^{\varDelta}_-(x,t)=&&\tilde{w}^{\varDelta}_-
		+\int^{x_{2N_x-1}}_{0}
		\zeta(u^{\varDelta}_{n,0}(x))dx+\int^x_{x_{2N_x-1}}
\zeta(\check{u}^{\varDelta}_-(y,t))dy\\&&&
		+\left\{g_2(x,t;\check{u}^{\varDelta}_-)+V(u_-)\right\}(t-t_n).
	\end{alignedat}\label{appro1}
\end{align}

\begin{remark}${}$
\begin{enumerate}
\item We notice that approximate solutions ${z}^{\varDelta}_-,{w}^{\varDelta}_-$ 
and $\tilde{z}^{\varDelta}_-,\tilde{w}^{\varDelta}_-$ correspond to 
	$z,w$ and $\tilde{z},\tilde{w}$ in \eqref{transformation}, respectively.
\item 
For $t_n<t<t_{n+1}$, our approximate solutions will satisfy 
\begin{align}
	\begin{alignedat}{2} \int^{x_{2N_x-1}}_{0}&
		\zeta(u^{\varDelta}(x,t_{n+1-}))dx+\int^{t_{n+1}}_{t_n}\sum_{0\leq x\leq x_{2N_x-1}}(\sigma[\eta_{\ast}]-[q_{\ast}])dt\\
		&=\int^{x_{2N_x-1}}_{0}
		\zeta(u^{\varDelta}_{n,0}(x))dx+V(u_-){\varDelta}t
		+o({\varDelta}x).
	\end{alignedat}
	\label{mass-conservation}	
\end{align}
\item Our construction of approximate solutions uses the iteration method twice (see \eqref{appro1-2} and \eqref{appro1}) 
to deduce \eqref{iteration}. 
\end{enumerate}

\end{remark}

We first set 
\begin{align*}
	\begin{alignedat}{2}
		\tilde{z}^{\varDelta}_+
				=z_+- \int^{x_{2N_x}}_{0}
		\zeta(u^{\varDelta}_{n,0}(x))dx,\;
		\tilde{w}_+
		=w_+- \int^{x_{2N_x}}_{0}
		\zeta(u^{\varDelta}_{n,0}(x))dx,
	\end{alignedat}
\end{align*}
where $z_+=-\dfrac{\left(\rho_+\right)^{\theta}}{\theta},\;w_+=\dfrac{\left(\rho_+\right)^{\theta}}{\theta}$.

We next construct $\check{u}^{\varDelta}_+$
\begin{align*}
\begin{alignedat}{2}
	\check{z}^{\varDelta}_+(x,t)&&=&\tilde{z}^{\varDelta}_+
+ \int^{x_{2N_x}}_{0}
		\zeta(u^{\varDelta}_{n,0}(x))dx
	+\int^x_{x_{2N_x}}\zeta(u_+)dy
+g_1(x,t;u_+)(t-t_n)
	,\\
	\check{w}^{\varDelta}_+(x,t)&&=&\tilde{w}^{\varDelta}_+
	+\int^{x_{2N_x}}_{0}
		\zeta(u^{\varDelta}_{n,0}(x))dx+\int^x_{x_{2N_x}}\zeta(u_+)dy
+g_2(x,t;u_+)(t-t_n).
\end{alignedat}
\end{align*}

Using $\check{u}^{\varDelta}_+(x,t)$, we define ${u}^{\varDelta}_+(x,t)$ as follows. \begin{align}
	\begin{alignedat}{2}
		{z}^{\varDelta}_+(x,t)&=&&\tilde{z}^{\varDelta}_+
+ \int^{x_{2N_x}}_{0}
		\zeta(u^{\varDelta}_{n,0}(x))dx+\int^x_{x_{2N_x}}\zeta(\check{u}_+(y,t))dy
		+g_1(x,t;\check{u}_+)(t-t_n)
		,\\
		{w}^{\varDelta}_+(x,t)&=&&\tilde{w}^{\varDelta}_+
		+ \int^{x_{2N_x}}_{0}
		\zeta(u^{\varDelta}_{n,0}(x))dx+\int^x_{x_{2N_x}}\zeta(\check{u}_+(y,t))dy
		+g_2(x,t;\check{u}_+)(t-t_n).
	\end{alignedat}
\label{appr-R}
\end{align}

\begin{figure}[htbp]
	\begin{center}
		\hspace{-2ex}
		\includegraphics[scale=0.3]{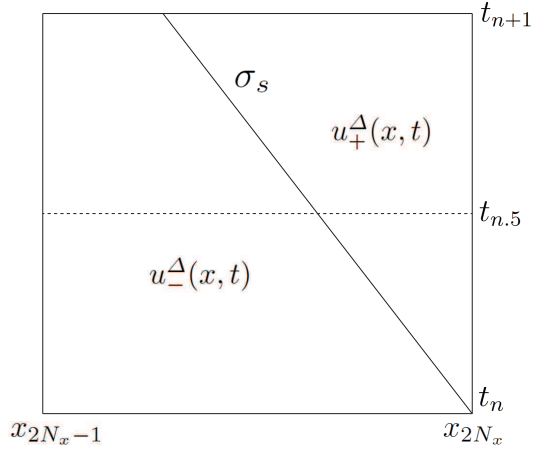}
	\end{center}
	\caption{Case 1: The case where a 1-shock wave arises near the boundary.}
	\label{Fig:case1}
\end{figure}

Then, we define approximate solution ${u}^{\varDelta}(x,t)$ 
in $D$ as follows (see Figure \ref{Fig:case1}).
\begin{align*}{u}^{\varDelta}(x,t)=
\begin{cases}
{u}^{\varDelta}_-(x,t),&D\cap\left\{x-x_{2N_x}\leq\sigma_s(t-t_n)\right\},\\ 
{u}^{\varDelta}_+(x,t),&D\cap\left\{x-x_{2N_x}>\sigma_s(t-t_n)\right\}.
\end{cases}
\end{align*}

\subsection{Case 2: the case where a rarefaction wave arises} Let $\alpha$ be a constant satisfying $1/2<\alpha<1$. Then we can choose 
a positive value $\beta$ small enough such that $\beta<\alpha$, $1/2+\beta/2<\alpha<
1-2\beta$, $\beta<2/(\gamma+5)$ and $(9-3\gamma)\beta/2<\alpha$.

{\it Step 1}.\\
In order to approximate a 1-rarefaction wave by a piecewise 
constant {\it rarefaction fan}, we introduce the integer  
\begin{align*}
	p:=\max\left\{[\hspace{-1.2pt}[(z_{+}-z_-)/({\varDelta}x)^{\alpha}]
	\hspace{-1pt}]+1,2\right\},
\end{align*}
where $z_-=z(u_-),z_{+}=z(u_{+})$ and $[\hspace{-1.2pt}[x]\hspace{-1pt}]$ is the greatest integer 
not greater than $x$. Notice that
\begin{align}
	p=O(({\varDelta}x)^{-\alpha}).
	\label{order-p}
\end{align}
Define \begin{align*}
	z_1^*:=z_-,\;z_p^*:=z_{+},\;w_i^*:=w_-\;(i=1,\ldots,p),
\end{align*}
and 
\begin{align*}
	z_i^*:=z_{2N_x-1}+(i-1)({\varDelta}x)^{\alpha}\;(i=1,\ldots,p-1).
\end{align*}
We next introduce the rays $x=1+\lambda_1(z_i^*,z_{i+1}^*,w_-)
(t-n{\varDelta}{t})$ separating finite constant states 
$(z_i^*,w_i^*)\;(i=1,\ldots,p)$, 
where  
\begin{align*}
	\lambda_1(z_i^*,z_{i+1}^*,w_-):=v(z_i^*,w_-)
	-S(\rho(z_{i+1}^*,w_-),\rho(z_i^*,w_-)),
\end{align*}
\begin{align*}
	\rho_i^*:=\rho(z_i^*,w_-):=\left(\frac{\theta(w_--z_i^*)}2\right)^{1/\theta}\;,
	\quad{v}_i^*:={v}(z_i^*,w_-):=\frac{w_-+z_i^*}2
\end{align*}
and

\begin{align}
	S(\rho,\rho_0):=\left\{\begin{array}{lll}
		\sqrt{\displaystyle{\frac{\rho(p(\rho)-p(\rho_0))}{\rho_0(\rho-\rho_0)}}}
		,\quad\mbox{if}\;\rho\ne\rho_0,\\
		\sqrt{p'(\rho_0)},\quad\mbox{if}\;\rho=\rho_0.
	\end{array}\right.
	\label{s(,)}
\end{align}

We call this approximated 1-rarefaction wave a {\it 1-rarefaction fan}.

\vspace*{10pt}
{\it Step 2}.\\
In this step, we replace the above constant states  with functions of $x$ and $t$ as follows:

In view of \eqref{transformation}, we construct ${u}^{\varDelta}_1(x,t)$. We first determine the approximation of $\tilde{z},\tilde{w}$ in \eqref{transformation} 
as follows.
\begin{align*}
	\begin{alignedat}{2}
		\tilde{z}^{\varDelta}_1
		=&z_--
\int^{x_{2N_x-1}}_{0}
		\zeta(u^{\varDelta}_{n,0}(x))dx,\;
		\tilde{w}^{\varDelta}_1=w_--
 \int^{x_{2N_x-1}}_{0}
		\zeta(u^{\varDelta}_{n,0}(x))dx.
	\end{alignedat}
	\end{align*}
We set       
\begin{align}
	\begin{alignedat}{2}
		&\check{z}^{\varDelta}_1(x,t)=&&\tilde{z}^{\varDelta}_1
		+\int^{x_{2N_x-1}}_{0}
		\zeta(u^{\varDelta}_{n,0}(x))dx
        +\int^x_{x_{2N_x-1}}\zeta(u_-)dy
\\&&&	
+\left\{g_1(x,t;u_-)+V(u_-)\right\}(t-t_n)
		,\\
		&\check{w}^{\varDelta}_1(x,t)=&&\tilde{w}^{\varDelta}_1
		+ \int^{x_{2N_x-1}}_{0}
		\zeta(u^{\varDelta}_{n,0}(x))dx+\int^x_{x_{2N_x-1}}\zeta(u_-)dy
\\&&&
+\left\{g_2(x,t;u_-)+V(u_-)\right\}
		(t-t_n).
	\end{alignedat}\label{appro1-2}
\end{align}
Using $\check{u}^{\varDelta}_1(x,t)$, we next define ${u}^{\varDelta}_1(x,t)$ as follows. 
\begin{align}
	\begin{alignedat}{2}
		&{z}^{\varDelta}_1(x,t)=&&\tilde{z}^{\varDelta}_1
+\int^{x_{2N_x-1}}_{0}
		\zeta(u^{\varDelta}_{n,0}(x))dx+\int^x_{x_{2N_x-1}}
\zeta(\check{u}^{\varDelta}_1(y,t))dy\\&&&	+
\left\{g_1(x,t;\check{u}^{\varDelta}_1)+V(u_-)\right\}(t-t_n),\\
		&{w}^{\varDelta}_1(x,t)=&&\tilde{w}^{\varDelta}_1
		+\int^{x_{2N_x-1}}_{0}
		\zeta(u^{\varDelta}_{n,0}(x))dx+\int^x_{x_{2N_x-1}}
\zeta(\check{u}^{\varDelta}_1(y,t))dy\\&&&
		+\left\{g_2(x,t;\check{u}^{\varDelta}_1)+V(u_-)\right\}(t-t_n).
	\end{alignedat}\label{appro1}
\end{align}

\begin{figure}[htbp]
	\begin{center}
		\hspace{-2ex}
		\includegraphics[scale=0.3]{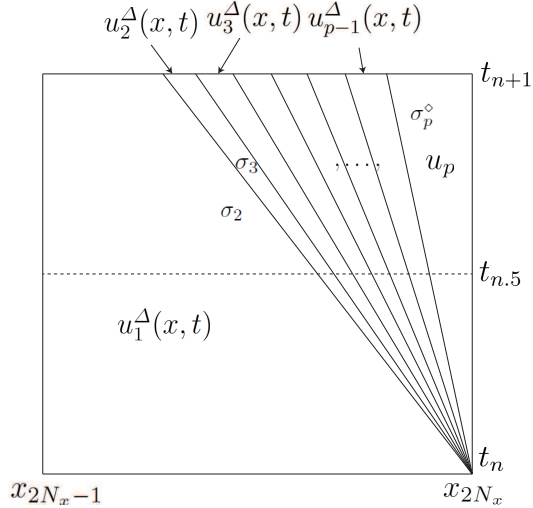}
	\end{center}
	\caption{Case 2: The case where a 1-rarefaction  arises near the boundary.}
	\label{Fig:case2}
\end{figure}

First, by the implicit function theorem, we determine a propagation speed $\sigma_2$ and $u_2=(\rho_2,m_2)$ such that 
\begin{itemize}
	\item[(1.a)] $z_2:=z(u_2)=z^*_2$
	\item[(1.b)] the speed $\sigma_2$, the left state ${u}^{\varDelta}_1(x^{\varDelta}_2(t_{n.5}),t_{n.5})$ and the right state $u_2$ satisfy the Rankine--Hugoniot conditions, i.e.,
	\begin{align*}
		f(u_2)-f({u}^{\varDelta}_1(x^{\varDelta}_2(t_{n.5}),t_{n.5}))=\sigma_2(u_2-{u}^{\varDelta}_1(x^{\varDelta}_2(t_{n.5}),t_{n.5})),
	\end{align*}
\end{itemize}
where $x^{\varDelta}_2(t)
=1+
\sigma_2(t-t_n)$. Then we fill up by ${u}^{\varDelta}_1(x)$ the sector where $t_n\leq{t}<t_{n+1},x_{2N_x-1}\leq{x}<x^{\varDelta}_2(t)$ (see Figure \ref{Fig:case2})
.

Assume that $u_k$, ${u}^{\varDelta}_k(x,t)$, a propagation speed $\sigma_k$ and $x^{\varDelta}_{k}(t)$  are defined. Then we similarly determine
$\sigma_{k+1}$ and $u_{k+1}=(\rho_{k+1},m_{k+1})$ such that 
\begin{itemize}
	\item[($k$.a)] $z_{k+1}:=z(u_{k+1})=z^*_{k+1}$,
	\item[($k$.b)] $\sigma_{k}<\sigma_{k+1}$,
	\item[($k$.c)] the speed 
	$\sigma_{k+1}$, 
	the left state ${u}^{\varDelta}_k(x^{\varDelta}_{k+1}(t_{n.5}),t_{n.5})$ and the right state $u_{k+1}$ satisfy 
	the Rankine--Hugoniot conditions, 
\end{itemize}
where $x^{\varDelta}_{k+1}(t)=1+\sigma_{k+1}(t-t_n)$. Then we fill up by ${u}^{\varDelta}_k(x,t)$ the sector where $t_n\leq{t}<t_{n+1},x^{\varDelta}_{k}(t)\leq{x}<x^{\varDelta}_{k+1}(t)$ 
.

We construct ${u}^{\varDelta}_{k+1}(x,t)$ as follows.

We first determine
\begin{align*}&
	\begin{alignedat}{2}
		&\tilde{z}^{\varDelta}_{k+1}=&&z_{k+1}-\int^{x_{2N_x-1}}_{0}
		\zeta(u^{\varDelta}_{n,0}(x))dx-V(u_-)\frac{{\varDelta}t}{2}
 -\sum^{k}_{l=1}\int^{x^{\varDelta}_{l+1}(t_{n.5})}_{x^{\varDelta}_l(t_{n.5})}
		\zeta(u^{\varDelta}_l(x,t_{n.5}))dx
,
\end{alignedat}\\
	&\begin{alignedat}{2}
&\tilde{w}^{\varDelta}_{k+1}=&&w_{k+1}-\int^{x_{2N_x-1}}_{0}
\zeta(u^{\varDelta}_{n,0}(x))dx-V(u_-)\frac{{\varDelta}t}{2}
-\sum^{k}_{l=1}\int^{x^{\varDelta}_{l+1}(t_{n.5})}_{x^{\varDelta}_l(t_{n.5})}
		\zeta(u^{\varDelta}_l(x,t_{n.5}))dx
,
	\end{alignedat}
\end{align*}
where $x^{\varDelta}_1(t)=x_{2N_x-1},\;x^{\varDelta}_l(t)=1+\sigma_l(t-t_n)\quad
(l=2,3,\ldots,k+1)$ and $t_{n.5}$ is defined in \eqref{terminology}.

We next define $\check{u}^{\varDelta}_{k+1}$ as follows.
\begin{align*}
	\begin{alignedat}{2}
		&\check{z}^{\varDelta}_{k+1}(x,t)=&&\tilde{z}^{\varDelta}_{k+1}+\int^{x_{2N_x-1}}_{0}
		\zeta(u^{\varDelta}_{n,0}(x))dx+V(u_-)(t-t_n)
		\\&&&+\sum^{k}_{l=1}
		\int^{x^{\varDelta}_{l+1}(t)}_{x^{\varDelta}_l(t)}
		\zeta(u^{\varDelta}_l(x,t))dx+\int^x_{x^{\varDelta}_{k+1}(t)}
		\zeta(u_{k+1})dy\\&&&
		+g_1(x,t;u_{k+1})(t-t_{n.5})
		+\int^{t}_{t_{n.5}}\hspace*{0ex}\sum_{\substack{x_{2N_x-1}\leq y \leq x}}\hspace*{-1ex}(\sigma[\eta_{\ast}]-[q_{\ast}])ds
		,\end{alignedat}
	\end{align*}\begin{align*}
	\begin{alignedat}{2}
		&\check{w}^{\varDelta}_{k+1}(x,t)=&&\tilde{w}^{\varDelta}_{k+1}
		+ \int^{x_{2N_x-1}}_{0}
		\zeta(u^{\varDelta}_{n,0}(x))dx+V(u_-)(t-t_n)
		\\&&&+\sum^{k}_{l=1}
		\int^{x^{\varDelta}_{l+1}(t)}_{x^{\varDelta}_l(t)}
		\zeta(u^{\varDelta}_l(x,t))dx
		+\int^x_{x^{\varDelta}_{k+1}(t)}
		\zeta(u_{k+1})dy\\&&&
		+g_2(x,t;u_{k+1})(t-t_{n.5})
		+\int^{t}_{t_{n.5}}\hspace*{0ex}\sum_{\substack{x_{2N_x-1}\leq y \leq x}}\hspace*{-1ex}(\sigma[\eta_{\ast}]-[q_{\ast}])ds.
	\end{alignedat}
\end{align*}

Finally, using $\check{u}^{\varDelta}_{k+1}(x,t)$, we define ${u}^{\varDelta}_{k+1}(x,t)$ as follows.

\begin{align}
	\begin{alignedat}{2}
		&{z}^{\varDelta}_{k+1}(x,t)&=&\tilde{z}^{\varDelta}_{k+1}+ \int^{x_{2N_x-1}}_{0}
		\zeta(u^{\varDelta}_{n,0}(x))dx+V(u_-)(t-t_n)\\&&&
		+\sum^{k}_{l=1}
		\int^{x^{\varDelta}_{l+1}(t)}_{x^{\varDelta}_l(t)}
		\zeta(u^{\varDelta}_l(x,t))dx
		+\int^x_{x^{\varDelta}_{k+1}(t)}
		\zeta(\check{u}^{\varDelta}_{k+1}(y,t))dy\\&&&
		+g_1(x,t;\check{u}^{\varDelta}_{k+1})(t-t_{n.5})
		+\int^{t}_{t_{n.5}}\hspace*{0ex}\sum_{\substack{x_{2N_x-1}\leq y \leq x}}\hspace*{-1ex}(\sigma[\eta_{\ast}]-[q_{\ast}])ds
		,\\
		&{w}^{\varDelta}_{k+1}(x,t)&=&\tilde{w}^{\varDelta}_{k+1}
		+ \int^{x_{2N_x-1}}_{0}
		\zeta(u^{\varDelta}_{n,0}(x))dx+V(u_-)(t-t_n)
		\\&&&+\sum^{k}_{l=1}
		\int^{x^{\varDelta}_{l+1}(t)}_{x^{\varDelta}_l(t)}
		\zeta(u^{\varDelta}_l(x,t))dx+\int^x_{x^{\varDelta}_{k+1}(t)}\zeta(
		\check{u}^{\varDelta}_{k+1}(y,t))dy\\&&&
		+g_2(x,t;\check{u}^{\varDelta}_{k+1})(t-t_{n.5})
		+\int^{t}_{t_{n.5}}\hspace*{0ex}\sum_{\substack{x_{2N_x-1}\leq y \leq x}}\hspace*{-1ex}(\sigma[\eta_{\ast}]-[q_{\ast}])ds.
	\end{alignedat}
\label{appr-k}
\end{align}

By induction, we define $u_i$, ${u}^{\varDelta}_i(x,t)$ and $\sigma_i$ $(i=1,\ldots,p-1)$.
Finally, we determine a propagation speed $\sigma_p$ and $u_p=(\rho_p,m_p)$ such that
\begin{itemize}
	\item[($p$.a)] $z_p:=z(u_p)=z^*_p$,
	\item[($p$.b)] the speed $\sigma_p$, 
	and the left state ${u}^{\varDelta}_{p-1}(x^{\varDelta}_{p}(t_{n.5}),t_{n.5})$ and the right state $u_p$ satisfy the Rankine--Hugoniot conditions, 
\end{itemize}where $x^{\varDelta}_{p}(t)=1+\sigma_{p}(t-t_n)$. 
We then fill up by ${u}^{\varDelta}_{p-1}(x,t)$ and $u_p$ the sector where
$t_n\leq{t}<t_{n+1},x^{\varDelta}_{p-1}(t)
\leq{x}<x^{\varDelta}_{p}(t)$ 
and the line $t_n\leq{t}<t_{n+1},x=x^{\varDelta}_{p}(t)$, respectively.

Given $u_{-}$ and $z_{+}$ with $z_{-}\leq{z}_{+}$, we denote 
this piecewise functions of $x$ and $t$ 1-rarefaction wave by 
$R_1^{\varDelta}(u_{-},z_{+},x,t)$.

 Finally, we construct ${u}^{\varDelta}_{p}(x,t)$ with 
 ${u}^{\varDelta}_{p}(1,t_n)=u_+$ 
in the similar manner to ${u}^{\varDelta}_+(x,t)$ in the 
Case 1 and fill up by ${u}^{\varDelta}_{p}(x,t)$ the sector where
$t_n\leq{t}<t_{n+1},x^{\varDelta}_{p}(t)
\leq{x}\leq1$.

\subsection{Case 3: the case where a rarefaction wave and the vacuum arise}

In this case, we consider the case where $\rho_{+}\leq({ \varDelta}x)^{\beta}$,
which means that $u_{+}$ is near the vacuum. In this case, we cannot construct 
approximate solutions in a similar fashion to the case 1--2. Therefore, 
we must
define $u^{ \varDelta}(x,t)$ in the different way.

\vspace*{5pt}
{\bf Case 3.1}
$\rho_->({ \varDelta}x)^{\beta}$\\
Let $u^{(1)}_-$ be a state satisfying $ w(u_-^{(1)})=w(u_-)$ and 
$\rho^{(1)}_-=({ \varDelta}x)^{\beta}$.

(i) $z(u_{+})-z(u^{(1)}_-)\leq({ \varDelta}x)^{\alpha}$\\
Notice that $w(u_+)=w(u_-)=w(u^{(1)}_-)$. 
Then there exists $C>0$ such that $\rho^{(1)}_--
\rho_+\leq{C}({ \varDelta}x)^{\alpha}$. Since $\alpha>\beta$, we then have 
$\rho_+\geq3({ \varDelta}x)^{\beta}/4$.
This case is reduced to the case 2.

(ii) $z(u_{+})-z(u^{(1)}_-)>({ \varDelta}x)^{\alpha}$\\
Set 
\begin{align*}
\bar{z}:=-M_{n+1}-L_n+ \int^{x_{2N_x-1}}_{0}
		\zeta(u^{\varDelta}_{n,0}(x))dx+V(u_-){\varDelta}t
		+\int^{x_{2N_x}}_{x_{2N_x-1}}\left\{\eta_{\ast}(u_-)+K\right\}dx.
\end{align*}
Let $u^{(2)}_-$ be a state connected to $u_-$ on the right by 
$R_1^{ \varDelta}(\max\{z^{(1)}_-,\bar{z}\})(u_-)$. Connecting the left 
and right states $u^{(2)}_-$, $u_+$ with $\rho_+=v_+=0$ by a rarefaction curve and the vacuum, we construct a Rimann solution $(u^{(2)}_-,u_+)$. 
Then, in the region where ${u}^{ 
\varDelta}(x,t)$ is $R_1^{ \varDelta}(\max\{z^{(1)}_-,\bar{z}\})(u_-)$, 
the definition of 
$u^{ \varDelta}(x,t)$ is similar to Case 2. In the other region, we 
define $u^{ \varDelta}(x,t)$ by the Riemann solution $(u^{(2)}_-,u_+)$ itself.

\vspace*{10pt}
{\bf Case 3.2} $\rho_-\leq({ \varDelta}x)^{\beta}$

\vspace*{5pt}
(i) $z(u_-)\geq{\bar{z}}$\\
In this case, we define $u^{ \varDelta}(x,t)$ as a Riemann solution 
$(u_-,u_+)$.

\vspace*{5pt}
(ii) $z(u_-)<\bar{z}$\\

Set 
\begin{align*}
	\bar{w}:=M_{n+1}+ L_n+\int^{x_{2N_x-1}}_{0}
	\zeta(u^{\varDelta}_{n,0}(x))dx+V(u_-){\varDelta}t
	-\int^{x_{2N_x}}_{x_{2N_x-1}}\nu \rho_-dx.
\end{align*}

Let $\lambda_1(u_-)$ be the 
1-characteristic speed of $u_-$.
In the region where 
$t_n\leq{t}<t_{n+1}$ and $x_{2N_x-1}
\leq{x}\leq 1+\lambda_1(u_-)(t-t_n)$,
we define $\bar{u}^{ \varDelta}(x,t)$ in the similar manner to 
$\bar{u}^{ \varDelta}_1(x,t)$ in Case 2.

We next take $u^{(3)}_-$ such that $z(u^{(3)}_-)=\max\{z_-,\bar{z}\}$ and 
$w(u^{(3)}_-)=\min\{w_-,\bar{w}\}$. We then solve a Riemann problem
$(u^{(3)}_-,u_+)$. In the region where 
$t_n\leq{t}<t_{n+1}$ and $1+\lambda_1(u_-)(t-t_n)<x\leq x_{2N_x}$, 
we define $\bar{u}^{ \varDelta}(x,t)$ as this Riemann solution.

\begin{remark}\label{rem:approximate}\normalfont
	The approximate solution $u^{\varDelta}(x,t)$ is piecewise smooth in each of the 
	divided parts of the cell. Then, in the divided part, $u^{\varDelta}(x,t)$ satisfies
	\begin{align*}
		(u^{\varDelta})_t+f(u^{\varDelta})_x=o(1).
	\end{align*}
\end{remark}

\section{The $L^{\infty}$ estimate of the approximate solutions}\label{sec:bound}

We deduce from (\ref{remark2.1}) the following
theorem:
\begin{theorem}\label{thm:bound}
	For $x_{2N_x-1}\leq x\leq 1$,
	\begin{align}
		\begin{alignedat}{2}
			&\displaystyle {z}^{\varDelta}(x,t_{n+1-})&\geq&-M_{n+1}-L_n
			+\int^x_{0}\zeta({u}^{\varDelta}(y,t_{n+1-}))dy-{\it o}({\varDelta}x),\\
			&\displaystyle {w}^{\varDelta}(x,t_{n+1-})
			&\leq& M_{n+1}+L_n+\int^x_{0}\zeta({u}^{\varDelta}(y,t_{n+1-}))dy+\int^{t_{n+1}}_{t_n}\sum_{0\leq x\leq1}(\sigma[\eta_{\ast}]-[q_{\ast}])dt\\&&&+{\it o}({\varDelta}x),
		\end{alignedat}
		\label{goal}
	\end{align}
	where 
$M_{n+1}$ is defined in \eqref{M_n}, 
$t_{n+1-}=(n+1){\varDelta}t-0$ and ${\it o}({\varDelta}x)$ depends only on the bound of solutions. 
\end{theorem}

\vspace*{5pt}

Now, in the previous section, we have constructed 
$u^{\varDelta}(x,t)$ near the boundary $x=1$. In this case, we are devoted to case 1 in particular. For case 2 and 3, we refer to \cite{T9} and \cite{T10}.

\subsection{Estimates of ${z}^{\varDelta}(x,t)$ for the case 
where a shock arises near the boundary}

We first consider $\tilde{z}^{\varDelta}_-$. We recall that
\begin{align*}
	\begin{alignedat}{2}
		\tilde{z}^{\varDelta}_-=z_-- \int^{x_{2N_x-1}}_{0}
		\zeta(u^{\varDelta}_{n,0}(x))dx.
	\end{alignedat}
\end{align*}
From \eqref{remark2.1}, we have $\tilde{z}^{\varDelta}_-\geq -M_n-L_n$.

Since 
\begin{align}
\check{u}^{\varDelta}_-(x,t)={u}^{\varDelta}_-(x,t)+
O(({\varDelta}x)^2),
\label{iteration}
\end{align}
recalling \eqref{mass-conservation}, we have 
\begin{align}
	\begin{alignedat}{2}
		&{z}^{\varDelta}_-(x,t)&=&\tilde{z}^{\varDelta}_-
		+ \int^{x_{2N_x-1}}_{0}
		\zeta(u^{\varDelta}_{n,0}(x))dx+V(u_-)(t-t_n)+\int^x_{x_{2N_x-1}}
\zeta(\check{u}^{\varDelta}_-(y,t))dy\\&&&
		+g_1(x,t;\check{u}^{\varDelta}_-)(t-t_{n})
	\\
		&&\geq&-M_n-L_n+ \int^{x_{2N_x-1}}_{0}
		\zeta(u^{\varDelta}_{n,0}(x))dx+V(u_-)(t-t_n)\\&&&+\int^x_{x_{2N_x-1}}
		\zeta({u}^{\varDelta}_-(y,t))dy
		+g_1(x,t;{u}^{\varDelta}_-)(t-t_{n})
	-o({\varDelta}x).
	\end{alignedat}
\label{section3-2}
\end{align}
If ${z}^{\varDelta}_-(x,t_{n+1-0})>-M_n-L_n+I^n_{2N_x-1}-\sqrt{{\varDelta}x}$, 
from \eqref{mass-conservation} and $M_{n+1}=M_n+O({\varDelta}x)$, we obtain  $\eqref{goal}_2$. Otherwise, 
from the argument \eqref{estimate1}, 
regarding $M_0-\delta t$ in \eqref{estimate1} as $M_n+L_n$, 
we have $g_1(x,t;{u}^{\varDelta}_-)>\delta$. 
From \eqref{mass-conservation}, we conclude $\eqref{goal}_1$.

Next, we next consider ${z}^{\varDelta}_+$. 
We introduce the following lemma holds. 
\begin{lemma}
There exists a unique piecewise smooth entropy solution 
$(\rho(x,t),$ $m(x,t))$ 
containing the vacuum state $(\rho=0)$ on $D$ for the 
problem \eqref{Euler} and \eqref{eqn:RiemannIB}
satisfying 
\begin{align*}
z(u(x,t))\geq\min(-w(u_-),z(u_-)),\;
w(u(x,t))\leq\max(w(u_-),0),\;
\rho(x,t)\geq0.
\end{align*}
\end{lemma}

In this case, it follows from \eqref{remark2.1} and 
the above lemma  that
\begin{align}
z({u}_+)\geq  
\min\{-M_n-L_n-I^n_{2N_x-1},-M_n-L_n+I^n_{2N_x-1}\}.
\label{section3-1}
\end{align}
On the other hand, we have
\begin{align}
	I^n_{2N_x-1}=I^n_{2N_x}+O({\varDelta}x),
	\label{I_N}
\end{align}
where $I^n_{2N_x}=\int^{x_{2N_x}}_{0}
\zeta(u^{\varDelta}_{n,0}(x))dx$.

Moreover, our approximate solutions satisfies the conservation of mass:
\begin{align}
\int^1_0\rho^{\varDelta}(x,t_{n-})dx
=\int^1_0\rho_0(x)dx+o(1)
\label{conservation-mass}
\end{align}
and the energy inequality: 
\begin{align}
	\int^1_0\eta_{\ast}(u^{\varDelta}(x,t_{n-}))dx\leq
	\int^1_0\eta_{\ast}(u_0(x))dx+o(1).
\label{energy-inequality}
\end{align}

From \eqref{boundary-bound2}, \eqref{I_N}--\eqref{energy-inequality}, we obtain
\begin{align}
	I^n_{2N_x}<-\mu+O({\varDelta}x).
\label{I_2N}
\end{align}

It follows from \eqref{section3-1} that 
\begin{align*}
	z({u}_+)\geq -M_n-L_n+I^n_{2N_x}
\end{align*}
by choosing ${\varDelta}x$ small enough.
Then, we have 
\begin{align*}
	\begin{alignedat}{2}
		\tilde{z}^{\varDelta}_+
				=z_+- \int^{x_{2N_x}}_{0}
		\zeta(u^{\varDelta}_{n,0}(x))dx\geq -M_n-L_n.
	\end{alignedat}
\end{align*}
Therefore, since $\check{u}^{\varDelta}_+(x,t)=u^{\varDelta}_+(x,t)+O(({\varDelta}x)^2)$, we conclude that
\begin{align*}
	\begin{alignedat}{2}
		{z}^{\varDelta}_+(x,t)&=&&\tilde{z}^{\varDelta}_+
		+ \int^{x_{2N_x}}_{0}
		\zeta(u^{\varDelta}_{n,0}(x))dx+\int^x_{x_{2N_x}}\zeta(\check{u}^{\varDelta}_+(y,t))dy
		+g_1(x,t;\check{u}^{\varDelta}_+)(t-t_n)\\
&\geq&&\tilde{z}^{\varDelta}_+
+ \int^{x_{2N_x}}_{0}
\zeta(u^{\varDelta}_{n,0}(x))dx+\int^x_{x_{2N_x}}\zeta({u}^{\varDelta}_+(y,t))dy
+g_1(x,t;{u}^{\varDelta}_+)(t-t_n)		.
	\end{alignedat}
\end{align*}
In this case, we can obtain $\eqref{goal}_1$ in the similar manner to \eqref{section3-2}. We can similarly obtain $\eqref{goal}_2$.

Our approximate solutions satisfy the following propositions holds (these proofs are similar to \cite{T1}--\cite{T3},\;\cite{T9}, \cite{T10}.).
\begin{proposition}\label{pro:compact}
The measure sequence
\begin{align*}
\eta_{\ast}(u^{\varDelta})_t+q(u^{\varDelta})_x
\end{align*}
lies in a compact subset of $H_{\rm loc}^{-1}(\Omega)$ for all weak entropy 
pair $(\eta_{\ast},q)$, where $\Omega\subset[0,1]\times[0,1]$ is any bounded
and open set. 
\end{proposition}
\begin{proposition} 
Assume that the approximate solutions $u^{\varDelta}$ are bounded and satisfy Proposition \ref{pro:compact}. Then there is a convergent subsequence $u^{\varDelta_n}(x,t)$
in the approximate solutions $u^{\varDelta}(x,t)$ such that
\begin{equation*}      
u^{\varDelta_n}(x,t)\rightarrow u(x,t)
\hspace{2ex}
\text{\rm a.e.,\quad as\;\;}n\rightarrow \infty.
\end{equation*} 
The function $u(x,t)$ is a global entropy solution
of the Cauchy problem \eqref{IP}.
\end{proposition}


\end{document}